\input amstex

\def\bcc#1{\Bbb C^{#1}}

\def\ee#1{e_{#1}}

\def\hd{, \hdots ,}

\def\ii{| II_X|}

\def\na{n+a}
\def\ooo#1#2{\omega^{#1}_{#2}}
\def\oo#1{\omega^{#1}}
\def\ot{\!\otimes\!}

\def\pp#1{\Bbb P^{#1}}
\def\ppp{\Bbb P}

\def\ql#1#2{q_{{#1} {#2}}}

\def\ra{\rightarrow}
\def\rl#1#2#3{r_{{#1}{#2}{#3}}}
\def\rr#1#2#3#4{r^{#1}_{{#2} {#3}{#4}}}

\def\tim{\text{Image}\,}
\def\tdim{\text{dim}\,}
\def\tcodim{\text{codim}\,}
\def\tsingloc{\text{singloc}\,}

\def\tmod{\text{ mod }}
\def\tmin{\text{ min }}

\def\upperp{{}^{\perp}}

\def\ww{\wedge}

\documentstyle{amsppt}
\magnification = 1200
\hsize =15truecm
\hcorrection{.5truein}
\baselineskip =18truept
\vsize =22truecm
\NoBlackBoxes
\topmatter
\title Is a linear space contained in a variety? - On the number of
derivatives needed to tell.
\endtitle
\rightheadtext{linear spaces on varieties}
\author
  J.M. Landsberg
\endauthor

\date { December 24, 1996}\enddate
\address{ Laboratoire de Math\'ematiques,
  Universit\'e Paul Sabatier, UFR-MIG
  31062 Toulouse Cedex 4,
  FRANCE}\endaddress
\email {jml\@picard.ups-tlse.fr }
\endemail
\thanks {Supported by  NSF grant DMS-9303704.}
\endthanks
\keywords {moving frames, osculating hypersurfaces,
            projective differential geometry, second fundamental forms}
\endkeywords
\subjclass{ primary 14, secondary 53}\endsubjclass
\abstract{ Let $X^n\subset\bcc\na$ or $X^n\subset\pp\na$ be a patch
of an analytic submanifold of an affine or projective space, let
$x\in X$ be a general point, and let $L^k$ be a linear space 
of dimension $k$ osculating
to order $m$
 at  $x$. If $m$ is large enough, one expects $L$ to be contained
in $X$ and thus $X$ contains a linear space
of dimension $k$
through  almost every point.
 We show that $L\subset X$ 
in the following cases:   $k=1$ and $m=n+1$;
  $k=n-1$, $a\geq 2$, and $m=2$;  $n\geq 4$,
$k=n-2$ and $m=4$.  We prove these results 
by first deriving the   order of osculation that
generically implies containment 
and then showing that in these 
cases containment must occur.
 If $X$ is a patch of a projective variety, we address
 the question as to whether   $X$ can be a smooth variety. We show that
if there is a $\pp k$ through each point and   $\tcodim(X)<\frac
k{n-k}$ then $X$ cannot be a smooth variety.}
\endabstract

\endtopmatter

\document

How many derivatives does 
one need to take to determine if a 
variety or analytic manifold is built out of linear spaces? 
It was  known classically (see e.g. [?]) that
 to see if a surface is ruled,
one needs to take three derivatives at a general point.
We generalize  this result  to $n$ dimensions in theorem 3 below.

In another direction, Ran [R], proved that if a patch of a projective 
variety $X^n$ is such that at
each point there is a line osculating to order $n+1$, then the union
of the lines has dimension at most $n+1$. He proved this as part
of a more general statement on $n+2$ secants, but, as he stated,
the case where the points all come together is the most difficult
(in fact the case of disjoint points is classical). Theorem 3 improves
Ran's
result to show that the union of the lines is in fact
at most $n$ dimensional, that is, the lines all lie on $X$.

One motivation for this work comes from the author's work
with B. Ilic on dual varieties ([BL1],[BL2]), as such varieties have linear
spaces through each point and one would like to understand
how these linear spaces fit together  by taking as few derivatives
as possible. Another motivation is the author's work on developing a
differential-geometric understanding of varieties that fail to be
complete intersections, [L], where the complementary problem of determining
when a variety is contained in a linear space arises.

\smallpagebreak

It is a pleasure to thank Mark Green for bringing the classical result
to my attention and for useful suggestions regarding the exposition.

\medpagebreak

\noindent \S 1.{\it Statements of results.}

\smallpagebreak

Throughout this paper,
 $X^n\subset\Bbb C^{\na}$ or $X^n\subset\pp\na$
 is a patch of an analytic submanifold
of an affine or projective space   and
  $x\in X$ is a general point. In particular, $X$ could
be   the smooth points of a projective variety.

Say $\tdim Y\leq \tdim X$ and $x\in X$. We will say $Y$ 
{\it osculates to order $0$ at $x$}  if $x\in Y$,
$Y$ {\it osculates to order $1$ at $x$}  if $\tilde T_xY\subseteq \tilde
T_xX$,
where $\tilde T_x$ denotes the embedded tangent space, and
{\it $Y$  osculates to order $m$ at $x$}  if
for all $v\in T_xY\subseteq T_xX$,
$F_k(v\hd v)=0$ for all $2\leq k\leq m$,
where 
$F_k$ denotes the $(k-2)$-nd
variation of the affine or projective
second fundamental form of $X$ at $x$
(the notation is such that $F_2=II$ and in this paper
we use the convention that $F_k$ includes
the information of the $k$-th fundamental form). 
See ([L], 2.18) for definitions
and explanations.   

If $Y$ is a curve, then
$Y$ osculates to order $m$ at $x$ if
 $X\cap Y$ has multiplicity $m+1$ at $x$. In general, $Y$
osculates to order $m$ at $x$ if every curve $C$ in $Y$ passing
through $x$ has multiplicity $m+1$ at $x$ (assuming $x\in Y$).

If one can show that there is a linear space $L$ contained in $X$ passing
through $x$, then since $x$ is a general point, in fact through almost all
points $y$ of $X$ there must be a linear space containing $y$ and
contained in $X$. Informally, one says that $X$ is
\lq\lq built out of linear spaces\rq\rq .

\proclaim{Theorem 1 (Generic Theorem)} Let
 $X^n\subset\Bbb C^{\na}$ or $X^n\subset\pp\na$ be a patch
of an analytic submanifold of an affine or projective space, and let
$x\in X$ be a general point.
 If $k\leq \frac n2$ and
$a[\binom{k+m}m -(k+1)]>k(n-k)$
or $k\geq\frac n2$ and
$a[\binom{2k-n+m}m -(2k-n+1)]>(2k-n)(n-k)$,
 then
any linear space
$L^k$, of dimension $k$, osculating to order $m$ at $x$ will be contained
in $X$ if
$F_2\hd F_{m}$ are sufficiently generic (in a sense to be made
precise in the proof).
Conversely, if 
$F_2\hd F_{m}$ are sufficiently generic, and $k,m$ are as
above, then one cannot tell if $L\subset X$ from less information.
\endproclaim

The improvement in the expectation when $k\geq\frac n2$ is explained
by the following theorem:

\proclaim{Theorem 2}Let $X^n\subset\Bbb C^{\na}$ or $X^n\subset\pp\na$ be a
patch
of an analytic submanifold of an affine or projective space, and let
$x\in X$ be a general point.
If there exists a linear space
$L^k$, of dimension $k$,
osculating to order two at $x$, then there exists
a linear subspace $M^{2k-n}\subset L^k$, of dimension $2k-n$, such that
$M^{2k-n}\subset X$.
\endproclaim

Theorem 1 implies that  one expects a line osculating to order
$\frac{n+1}a-1$ to be contained in $X$.  
When $m=2$, the genericity condition
for a line is easy to state:
Let $v\in T_xX$ be a tangent vector in the direction of the line.
Let $II_v : T_xX/\{v\}\ra N_xX$
be the linear map
obtained by contracting the second
fundamental form with $v$.
The condition is that $II_v$
is injective. Note that this is indeed
a genericity condition when $a-1\geq n$ (and  not possible
for smaller $a$).  We prove that in any case, osculation to order
$m+1$ is enough:

\proclaim{Theorem 3}
Let $X^n\subset\Bbb C^{\na}$ or $X^n\subset\pp\na$ be a patch
of an analytic submanifold of an affine or projective space, and let
$x\in X$ be a general point.
 Say a line $L^1$ osculates to order $n+1$ at $x$.
Then $L\subset X$.
\endproclaim

In the other extreme, one expects an $L^{n-1}$ osculating to order
two to be contained in $X$. We prove

\proclaim{Theorem 4}Let $X^n\subset\Bbb C^{\na}$ or $X^n\subset\pp\na$ be a
patch
of an analytic submanifold of an affine or projective space, and let
$x\in X$ be a general point.
 If $a\geq 2$,
$X$ is not contained in a hyperplane, and a linear space $L^{n-1}$
osculates to order two at $x$, then $L\subset X$.
In particular, any variety with the second fundamental form
of a scroll is a scroll.
\endproclaim

Theorem 1 implies that in any codimension, if $k\geq \frac n2$, one expects
a $L^k$ osculating to order four to be contained in $X$. 
We show this expectation
holds for $n=4$, in fact, we show:

\proclaim{Theorem 5} Let $X^n\subset\Bbb C^{\na}$ or $X^n\subset\pp\na$
 be a patch
of an analytic submanifold of an affine or projective space, and let
$x\in X$ be a general point.
If $n\geq 4$ and a linear space
$L^{n-2}$  osculates to order four at $x$, then
$L^{n-2}\subset X$.
\endproclaim

Often when a projective variety contains large linear spaces, it is forced
to be singular.
In the case of a scroll, or more generally a variety
$X$ which is a fibration with linear fibers, or more generally
yet, the case $X^n\subset \pp\na$ is a fibration with any sort of fiber
of dimension $f$, then
in order for $X$ to be smooth
one must have  $f\leq a$ (see [RV]).
If the linear spaces overlap (e.g. on the quadric hypersurface),
  $X$ can still be smooth in small  codimension. We prove:

\proclaim{Theorem 6}
Let $X^n\subset\pp\na$ be a smooth variety
such that through each general point
  $x\in X$ there is at least one $\pp k\subset X$ containing $x$.
Then $a\geq k/(n-k)$.
\endproclaim
 
\bigpagebreak

\noindent \S 2. {\it Proofs}.

  In the case of the
results independent of codimension, there is no loss of generality
in assuming $X$ is a hypersurface, as projections to a hypersurface show
that the linear space osculates to infinite order to each tangent
hyperplane.

\demo{Proof of theorem 2}Assume $X$ is a hypersurface.
Thus its second fundamental form consists of a single quadric which
has at least a $k\geq \frac n2$ dimensional linear space in its base locus,
therefore its rank can be at most $2(n-k)$. 

But now the singular locus of $II$  corresponds directions
in the fiber of the Gauss mapping, which is   a linear space of
dimension at least $n-2(n-k)$.\qed
\enddemo

In what follows, anytime there is a singular locus in $II$,
the directions in the singular locus can safely be ignored, since
they already correspond to tangent directions to a linear space contained
in $X$ and an easy computation shows that if directions 
in the singular locus are
a subspace of a subspace tangent to an osculating linear space,
then if any complement to the singular directions osculate
to infinite order, the span of all the directions will. 

\demo{proof of  theorem 1}
Let
$F$ be the  subspace of the fiber of the Gauss map
contained in $L$ mandated by theorem 2,
and let $W:=T_xL/T_xF$. Let
$(v_1\hd v_n)=(v_{\xi},v_{\rho},v_{\psi})$ denote a basis
of $T_xX$ adapted such that
$T_xF= \{   v_{\psi} \}$
and  $W= \{ v_{\xi}\} \tmod \{ v_{\psi} \}$.
The index ranges are $1\leq \xi,\eta\leq \tmin\{k,n-k\}$,
$\tmin\{k+1,n-k+1\}\leq\rho,\sigma\leq\tmin\{n, 2(k-n)\}$,
$2(k-n)\leq \psi\leq n$ where the last range may be empty.
Following the remark above, we suppres the $\psi$ index
in all formulae.

By hypothesis $F_j(v_{\xi_1}\hd v_{\xi_j})=0$ for all $1\leq j\leq m$.
The coefficients of $F_{i+1}$ in this range are given by the formula
([L],2.20) which simplifies under our hypotheses to
$$
\rr m{\xi_1\hd}{\xi_i}\rho\oo\rho =
\frak S_{\xi_1\hd \xi_i}\rr m{\xi_1\hd}{\xi_{i-1}}\sigma\ooo\sigma{\xi_i}.
\tag 1
$$
Let $M=(T_xL)\upperp$. We have a series of maps
$$
R_j :
S^jW\ot N\ra W\ot M \tmod \tim R_{j-1}. \tag 2
$$
It is a genericity condition on 
$F_2\hd F_{m}$ that these maps all have maximal dimensional image
and this is the genericity condition refered to in the statement of
the theorem. Assuming this condition holds, we see
that   $\ooo\sigma\eta\equiv 0\tmod \{\oo\rho\}$ $\forall \sigma,\eta$.
By induction, for $i>m$ we have 
$$
\rr m{\xi_1\hd}{\xi_i}{\xi_{i+1}}\oo{\xi_{i+1}} +
\rr m{\xi_1\hd}{\xi_i}\rho\oo\rho =
\frak S_{\xi_1\hd \xi_i}\rr m{\xi_1\hd}{\xi_{i-1}}\sigma\ooo\sigma{\xi_i},
$$
which implies
$\rr m{\xi_1\hd}{\xi_i}{\xi_{i+1}}=0$
for all $i$, and thus the linear space
osculates to infinte order.

The first assertion of the theorem follows by noting that
$\tdim N=a$, $\tdim W= \tmin\{ k,2k-n\}$ and $\tdim M=  n-k $, and
the  converse follows  by noting that we have used all the terms in the
$F_j$ that are subject to restrictions.
\qed\enddemo

\demo{proof of theorem 3} Assume we are
in the case of codimension one.
In this case the map (2) is always from a one
dimensional vector space. We claim that if   $R_j=0$ for some $j$,
then all higher maps are zero and the subspace of $M$ missed never
appears in computing the coefficients.

 More precisely, if the maps in question are not all nonzero, 
  suppress the
normal index in $F_j$ but use the notation
$r^{(j)}$ to denote the coefficients of $F_j$. Introduce additional indices

$b,c$ and $x,y$ splitting the $\rho$ indices such that at orders up to
$n+1$,
$r^{(j)}_{1\hd 1b}=0$ and the $R_j$'s surject onto the space
with the $x$ indices. Thus
 $\ooo x1\equiv 0\tmod \{\oo y\}$. Then at order $n+2$,  
$$
r^{(n+2)}_{1\hd 11}\oo  1
+r^{(n+2)}_{1\hd 1b}\oo  b + r^{(n+2)}_{1\hd 1x}\oo  x
=r^{(n+1)}_{1\hd 1y}\ooo y1,
$$
which implies that $r^{(n+2)}_{1\hd 11}=r^{(n+2)}_{1\hd 1b}=0$
and one continues by induction.\qed\enddemo

\demo{proof of theorem 4}
Without loss of generality, assume $a=2$. 
By the proof of theorem 2, at most three variables appear in $\ii$.
Up to equivalence, there are four cases for
$II_{X,x}$. 

Case 1: One can normalize such that
$q^{n+1}= \oo 1\oo n, q^{n+2}=\oo n\oo n$.
In this case
$
\rr{n+2}11\beta \oo \beta =0
$
which implies
$$
\rr{n+2}1n\beta\oo\beta=\rr{n+2}1nn\oo n=\ooo n1
$$
i.e. that
$$
\ooo n1\equiv 0\tmod \{\oo n\}
$$
and we are done by the proof of theorem 1.

Case 2: One can normalize such that
$q^{n+1}= \oo 1\oo n, q^{n+2}=0$. This case is not possible,
as since $X$ is not contained in a hyperplane, one would need
 the third fundamental
form to be nonzero, which is impossible as the prolongation of this
system of quadrics is empty ([GH], 1.47).

Case 3: One can normalize such that
$q^{n+1}= \oo 1\oo n, q^{n+2}=\oo 2\oo n$. In this case
$$
\rr{n+1}22\beta =\rr{n+2}11\beta=0
$$
which implies
$$
\align
&\rr{n+1}121\oo 1 + \rr{n+1}12n\oo n =\ooo n2 \\
&\rr{n+2}122\oo 2 + \rr{n+2}12n\oo n =\ooo n1\endalign
$$
i.e. that $\ooo n1\equiv 0\tmod \{ \oo 2,\oo n \}$
and  $\ooo n2\equiv 0\tmod \{ \oo 1,\oo n \}$.
But now
$$
\align
&\rr{n+1}111\oo 1 + \rr{n+1}11n\oo n =2\ooo n1 \\
&\rr{n+2}222\oo 2 + \rr{n+2}22n\oo n =2\ooo n2\endalign
$$
which implies $\{\ooo n1,\ooo n2\}\equiv 0\tmod \{\oo n\}$.

Case 4: One can normalize such that $q^{n+1}= (\oo n)^2$,
$q^{n+2}=0$. In this case $L$ is the fiber of the Gauss map.

In all cases, to see $X$ is a scroll, i.e., the total space of a projective
bundle
$\ppp (E)$
over a curve $C$ with fibers embedded linearly, let
$C$ be an integral curve of the vector field $\ee n$ through $x$,
and let $\ppp E_x=L$.
\qed\enddemo

\demo{Proof of theorem 5}
  Use indices
 as above and
normalize $q$ such that
$
\ql\xi{(m+\eta)} = \delta_{\xi\eta}$ where
$m=\tmin \{ k, 2k-n\}$.
We are given that for all $\xi_1\hd\xi_4$,
$\rl{\xi_1}{\xi_2}{ \xi_3}=0,$
$\rl{\xi_1}{\xi_2}{\xi_3 \xi_4}=0$.
Thus
$$
\align
\rl\xi\eta\rho \oo\rho
&=\ooo{m+\xi}\eta +\ooo{m+\eta}\xi,\\
\rl{\xi_1}{\xi_2}{\xi_3\rho}\oo\rho
&=\frak S_{ \xi_1 \xi_2 \xi_3}\rl{\xi_1}{\xi_2}{(m+\eta)}
\ooo{m+\eta}{\xi_3}\endalign
$$
Let $W=\{v_1\hd v_m\}$  be as above and consider 
$\rl{\xi_1}{\xi_2}{(m+\eta)}$ as giving a system of quadrics $R$
on
$W$ parametrized by 
$q(W)=\{v_{m+1}\hd v_{2m}\}$ which we also identify
with $W^*$ via $q$. Thus  
$$
 R_3:  S^3W \ra q(W)^*\ot W \tmod \tim R_2
$$
  quotients to a map
$$
\align
\tilde R: S^3W &\ra \Lambda^2 W. \\
v_{\xi_i}v_{\xi_j}v_{\xi_k}&\mapsto
\frak S_{ijk}\rl {\xi_i}{\xi_j}{(m+\eta )}
v^{\eta}\ww v^{\xi_k}\endalign
$$

If $\tilde R$ is surjective, then we are done. 

Claim: if  $R$ is an $m$-dimensional system of quadrics,
 then $\tilde R$ is surjective. In fact
$\tilde R$ is surjective if $R$ is an
$(m-1)$-dimensional system. To see this note that if one writes
out the matrix representing $R_3$, the row corresponding
to $v_{\eta}\ot v^{n+\zeta}$ consists of the entries 
$\rl {\xi_1}{\xi_2}{n+\zeta}$
in the slot $(\xi_1,\xi_2,\eta )$, so it is clear that if the quadrics
are all independent, the map is surjective.  However, since we only
care about the map to $\Lambda^2W^*$, the map is still surjective if
there are only $(m-1)$ independent quadrics.

On the other hand, if the map is zero
then we are also done because one has
$$
\rl{\xi_1\hd\xi_g}\beta{}\oo\beta
=0
$$
again by induction. So in any case, if $k=n-2$, $L\subset X$.\qed\enddemo

\noindent{\bf Remark}: It is often the case, 
as suggested by the paragraph above, that   osculation
to the expected order can only fail to imply containment if a
condition between degeneracy and genericity holds.
 
\demo{proof of theorem 6}
In order for $X$ to be smooth, the singular locus of the second fundamental
form must be empty. The proof will follow from the following lemma:

\proclaim{Lemma} Let $A$ be an $a$-dimensional system of quadrics
  on an $n$-dimensional vector space $V$, with a linear space
$W$  of dimension $k$ in
the base locus of $A$. If
$a<k/(n-k)$, then $A$ has a singular locus.
\endproclaim

\demo{Proof}
Any quadric $q\in A$ can be written
$$
q= v^1w^1+\hdots + v^kw^k + q'
$$
where,   $q'\in S^2W\upperp$,
and $v^j\in W\upperp$. Since $k>n-k$, at most 
$n-k$ of the $v^j$ are independent.
Thus each quadric has at least an
$k-(n-k)=2k-n$ dimensional singular locus in $W$, so
if $a(n-k)<k$ then $\tsingloc (A)\neq 0$.\qed\enddemo\enddemo

\enddocument